\documentclass[11pt]{article}
\usepackage{amssymb}
\usepackage{bbm}
\newcommand{\Boxneu}{\square}
\newtheorem{pr}{Proposition} 
\newtheorem{de}{Definition}

\begin{document}
\title{Locally trivial quantum vector bundles and associated vector
bundles}
\author{Dirk Calow$^1$\thanks{supported by Deutsche Forschungsgemeinschaft,
e-mail
Dirk.Calow@itp.uni-leipzig.de}
\ and Rainer Matthes$^{1,2}$ \thanks{supported by S\"achsisches
Staatsministerium   f\"ur Wissenschaft und Kunst,\newline\hspace*{.5cm} e-mail
Rainer.Matthes@itp.uni-leipzig.de or rmatthes@mis.mpg.de}\\[.5cm]
\normalsize $^1$Institut f\"ur Theoretische Physik
\normalsize der Universit\"at Leipzig\\
\normalsize Augustusplatz 10/11,
\normalsize D-04109 Leipzig,
\normalsize Germany
\\[.5cm]
\normalsize $^2$Max-Planck-Institut f\"ur Mathematik
\normalsize in den Naturwissenschaften\\
\normalsize Inselstra{\ss}e 22-26,
\normalsize D-04103 Leipzig,
\normalsize Germany}
\date{}
\maketitle
\begin{abstract}
We define locally trivial quantum vector 
bundles (QVB) and construct such QVB associated to locally trivial quantum
principal 
fibre bundles. The 
construction is quite analogous to the classical construction of 
associated bundles. A covering of such bundles is induced from the
covering of the subalgebra of coinvariant elements of the principal bundle.
There exists a differential
structure on the associated vector bundle coming from the differential 
structure on the 
principal bundle, which allows to define connections on the associated
vector bundle associated to connections on the principal bundle.
\end{abstract}
This is the third in a series of papers devoted to locally trivial quantum 
bundles, following \cite{cama} and \cite{cama1}. 
The main aim of the
present paper is to define associated bundles and connections in the scheme of
locally trivial quantum principal bundles of \cite{Kon} and \cite{cama1}.
Again we follow the idea of gluing all objects from locally given objects.
We note that associated bundles have been defined in \cite{dur1}, \cite{dur1.1}
and \cite{hama1} as associated bimodules of colinear maps (intertwiners)
and in \cite{brma}, \cite{Kon}, \cite{pf1} as cotensor products. Both
definitions are equivalent (by some duality argument). We use the second
definition.

We start with some general considerations about covering und gluing of modules.
Then we define QVB over algebras with a complete covering. They have as typical
fibre a usual vector space, and possess local trivializations with suitable
properties. Any LC differential algebra on the basis gives rise to a certain
``differential structure'' on the QVB, which is an analogue of a module
of differential form valued sections in the classical case. Given such
a structure, we define the notions of connection (as covariant derivative) and
curvature on a QVB. Given a locally trivial quantum principal bundle and a
left coaction of its structure group (Hopf algebra) on some vector space, we
define an associated QVB as a cotensor product. This is indeed a locally
trivial QVB in the sense of our definition, whose transition functions come
from the transition functions of the principal bundle in the usual way.

As is known from \cite{cama1}, there is a maximal embeddable LC differential
algebra related to the differential structure of a locally trivial quantum
principal bundle. The differential structure on the QVB defined by this
LC differential algebra is isomorphic to the module of horizontal forms
``of type $\rho$'' on the principal bundle, where $\rho$ is the 
left coaction defining the associated bundle. 
Finally, we show that to every connection on the principal bundle there can
be associated a connection on the associated QVB. Locally, the associated
connection and its curvature are nicely determined by the local connection
and curvature forms of the principal connection.

\section{Coverings of modules}
First we give a definition of modules over an algebra being equivalent to
the
usual one. 
For a vector space $E$ we denote by $End(E)$ the algebra of linear
endomorphisms 
of $E$.
\begin{de} Let $E$ be a vector space, let $B$ be an algebra and $\kappa : B
\longrightarrow End(E) $ a linear map. \\$(E,B,\kappa)$ is called left
module if $\kappa$ 
satisfies $\kappa(ab)=\kappa(a) \kappa(b),\ ;a,b \in B$.
\\ $(E,B,\kappa)$ is called right module if $\kappa$ satisfies 
$\kappa(ab)=\kappa(b) \kappa(a)$.
\\For a linear subspace $Q \subset E$, $(Q,B,\kappa)$ is a submodule if 
$\kappa(B)(Q) \subset Q$. \end{de}
\begin{de} Let $(E,B,\kappa)$ be a left (right) module. $(E,B,\kappa)$ is
called faithful if $ker \kappa =0$. \end{de}
In analogy to the case of algebras one can 
define coverings of modules.
They will be needed in the definition of quantum vector bundles.
\begin{de} Let $(E,B,\kappa)$ be a left (right) module and let $(Q_i)_{i
\in I}$ be a finite family of left (right) submodules of $E$. $(Q_i)_{i \in
I}$ is called covering of $E$ if $\bigcap_i Q_i =0$. \end{de}
Obviously, as in the case of algebras, for a given covering $(Q_i)_{i \in
I}$ of a module $(E,B,\kappa)$ one obtains a family of vector spaces 
$E_i := E/Q_i$ with corresponding projections $q_i: E \longrightarrow 
E_i$. Since $Q_i$ are submodules, there exist linear maps 
$\kappa_i : B \longrightarrow End(E_i)$ defined by \[ \kappa_i(a) \circ q_i
=q_i 
\circ \kappa(a)\; a \in B \] such that $(E_i,B,\kappa_i)$ are left or right
modules respectively.
\\\\Remark: Since $\kappa_i$ are homomorphisms for left modules and 
antihomomorphism for right modules, $ker \kappa_i$ are ideals in $B$.
\begin{de} A covering $(Q_i)_{i \in I}$ of a left (right) module 
$(E,B,\kappa)$ is called nontrivial if $ ker \kappa_i \not= 0~\forall i \in
I$. \end{de}
\begin{pr} Let $(E,B,\kappa)$ be a faithful left (right) module and let 
$(Q_i)_{i \in I}$ be a covering of $E$. \\Then $(ker \kappa_i)_{i \in 
I}$ is a covering of $B$. 
\\The covering $(Q_i)_{i \in I}$ is nontrivial if the covering $(ker 
\kappa_i)_{i \in I}$ is nontrivial. \end{pr}
Proof: We have to prove that $\bigcap_i ker \kappa_i =0$.
It is easy to verify that \[ ker \kappa_i= \{ a \in B|\,\kappa(a)(E)
\subset Q_i).\] Now it is clear that \[\bigcap_i ker \kappa_i= \{ a \in
B|\,\kappa(a)(E) \subset Q_i \;\forall i \in I\} \] and since 
$\bigcap_i Q_i=0$
and $(E,B,\kappa)$ is faithful it follows $\bigcap
ker\kappa_i=0$.\hfill$\Boxneu$
\\\\\\Let $B_i:= B/ker \kappa_i$. It is obvious that the modules 
$(E_i,B_i,\tilde{\kappa_i})$, where $\tilde{\kappa_i}$ is defined by 
\[ \tilde{\kappa}_i \circ \pi_i = \kappa_i, \] are faithful.
Let \begin{eqnarray*} q_{ij}:E & \longrightarrow &  E/(Q_i+Q_j):=E_{ij}\\
q^i_j: E_i & \longrightarrow & E_{ij} \end{eqnarray*} be the canonical 
projections.
Assume that $(E,B,\kappa)$ is a faithful left (right) module and $(Q_i)_{i
\in I}$ is a covering of $E$. One has the vector space 
\begin{equation} E_c:=\{ (e_i)_{i \in I} \in \bigoplus_{i \in I} 
E_i|\; q^i_j(e_i)=q^j_i(e_j)\} \end{equation} and an injective 
homomorphism \[ {\cal K} : E \longrightarrow E_c \] by ${\cal 
K}(e)=(q_i(e))_{i \in I} $.
\begin{pr} Let $(E,B,\kappa)$ be a faithful left (right) module and let 
$(Q_i)_{i \in I}$ be a covering of $E$. Let $B_c$ be the covering
completion 
of $B$ with respect to $(ker \kappa_i)_{i \in I}$. Then there exists a
linear
map $\kappa_c : B_c \longrightarrow End(E_c)$ such that
$(E_c,B_c,\kappa_c)$ 
is a faithful left (right) module satisfying 
\[ {\cal K}(\kappa(a)(e)) = \kappa_c(K(a))({\cal K}(e)), \; a \in B, \;e
\in E,\] \end{pr} where $K:B \longrightarrow B_c$ is the injektive Homomorphism 
defined by $K(a)=(\pi_i(a))_{i \in I}$.
Proof: Since $(Q_i)_{i \in I}$ is a family of submodules, there 
exist linear maps $\kappa_{ij}: B \longrightarrow End(E_{ij}) $ defined by
\[ \kappa_{ij}(a) \circ q_{ij} = q_{ij} \circ \kappa(a),\; a \in B \] such
that $(E_{ij},B,\kappa_{ij})$ is a left (right) module. Let $B_{ij}:=
B/(ker 
\kappa_i + ker \kappa_j)$. Now one can define the linear map 
$\tilde{\kappa}_{ij} : 
B_{ij} \longrightarrow End(E_{ij}) $ by the formula 
\[ \tilde{\kappa}_{ij} \circ \pi_{ij} = \kappa_{ij}, \] thus
$(E_{ij},B_{ij}, \tilde{\kappa}_{ij})$ is a left (right) module. It is easy
to verify that $\tilde{\kappa}_{ij}$ satisfies \[ q^i_j \circ \tilde{\kappa}_i(a) 
=\tilde{\kappa}_{ij}(\pi^i_j(a)) \circ q^i_j,\;a \in B_i.\] Recall that
$B_c:=\{(a_i)_{i \in I} \in \bigoplus_{i \in I} B_i|
\pi^i_j(a_i)=\pi^j_i(a_j)\}$. 
We define $\kappa_c : B_c \longrightarrow End(E_c)$ as follows:
\begin{equation} \kappa_c((a_i)_{i \in I})((e_i)_{i \in 
I})=(\tilde{\kappa}_i(a_i)(e_i))_{i \in I}. \end{equation}
One has to prove that the image of $\kappa_c(a_c)$ lies in $E_c$ for all
$a_c \in B_c$: \begin{eqnarray*} 
q^i_j(\tilde{\kappa}_i(a_i)(e_i)) &=& 
\tilde{\kappa}_{ij}(\pi^i_j(a_i))(q^i_j(e_i)) \\
&=& \tilde{\kappa}_{ij}(\pi^j_i(a_j))(q^j_i(e_j)) \\
&=& q^j_i(\tilde{\kappa}_j(a_j)(e_j)). \end{eqnarray*}
The other properties of $\kappa_c$ follow from the definition. \hfill$\Boxneu$
\begin{de} Let $(E,B,\kappa)$ be a faithful left (right) module and let 
$(Q_i)_{i \in 
I}$ be a covering of $E$. The covering $(Q_i)_{i \in I}$ is called complete, if the 
family of ideals $(ker \kappa_i)_{i \in I}$ is a complete covering of $B$
and 
the injective linear map ${\cal K} :E \longrightarrow E_c$ is a left
(right) module isomorphism. 
\end{de}

\section{Locally trivial quantum vector bundles and associated vector
bundles}

On the algebraic level, the notion of vector bundle is in the classical
case 
related to the notion of section of a vector bundle. Let $M$ be a
manifold,
let $C(M)$ be the algebra of continous functions over $M$ and let $V$ be a 
vector space. The corresponding trivial vector bundle has the form $M
\times V$. 
It is known, that the set of sections of $M \times V$ is the set of all 
$V$-valued functions denoted by $C(M) \otimes V$. This classical background
leads us to the following definition of a locally trivial vector bundle.

\begin{de}
A locally trivial quantum vector bundle (QVB) is a tupel 
\begin{equation} \{(E,B,\kappa),V,(\zeta_i,J_i)_{i \in I}\} \end{equation} 
where $(E,B,\kappa)$ is a faithful left module, $(J_i)_{i \in I}$ is a 
complete covering of $B$, $V$ is a vector space and 
$\zeta_i : E \longrightarrow B_i \otimes V$ are surjective left module 
homomorphisms with 
the properties \begin{eqnarray}  (ker \zeta_i)_{i \in I}  \; 
&& \mbox{complete 
covering of E} \\  \zeta_i(\kappa(a)(e)) &=& \pi_i(a)\zeta_i(e) \;\;a \in
B\;e 
\in E \\ \label{kerzet} ker \zeta_i + ker \zeta_j &=& ker ((\pi^i_j \otimes
id) \circ \zeta_i)= ker((\pi^j_i \otimes id) \circ \zeta_j). 
\end{eqnarray} 
\end{de}
Remark: In this definition 
we have used the left module structure $(B_i \otimes 
V),B,\kappa_i)$, which is defined by \[\kappa_i(a)(b_i \otimes
v)=\pi_i(a)b_i 
\otimes v \;\; a \in B,\;b_i \in B_i,\; v \in V.\]
In the following we want to denote such a vector bundle by $E$.
\\By definition, for a locally trivial QVB $E$ the family of submodules
$(ker 
\zeta_i)_{i \in I}$ is a complete covering of $E$, i.e. $E$ is isomorphic
to 
its 
covering completion.
Note that there are isomorphisms $\tilde{\zeta}_i: E/ker \zeta_i
\longrightarrow 
B_i \otimes V$ 
defined by \[ \tilde{\zeta}_i(e + ker\zeta_i)=\zeta_i(e). \]
Because of (\ref{kerzet}) there exist also isomorphisms $\zeta^i_{ij} :
E/(ker 
\zeta_i + ker \zeta_j) \longrightarrow B_{ij} \otimes V$ defined by 
\[ \zeta^i_{ij}(e + ker \zeta_i + ker\zeta_j)=\zeta_i(e) + \zeta_i(ker
\zeta_j), 
\] thus there are left $B_{ij}$-module isomorphisms $\phi_{{ij}_E} $
defined by 
$\phi_{{ij}_E}:= 
\zeta^i_{ij} \circ {\zeta^j_{ij}}^{-1}$ such that the covering completion
of $E$ 
has the form
\[ E=\{ (e_i)_{\ i \in I} \in \bigoplus_{i \in I} B_i \otimes V| (\pi^i_j 
\otimes id)(e_i)=\phi_{{ij}_{E}} \circ (\pi^j_i \otimes id)(e_j) \}. \]
\begin{pr} Let $N$ be a locally trivial topological vector bundle over a
compact
topological 
space $M$.\\
Then the set 
of continuous sections $\Gamma(N)= \{ s: M \longrightarrow N\}$ is a
locally trivial
QVB. 
\\
Let $B=C(M)$ be the algebra of continuous functions over a compact
topological space
$M$. Let $(J_i)_{i\in I}$ be a finite covering of $B$ coming from a finite
covering
$(U_i)_{i\in I}$ of $M$ by closed sets with nonempty open interior. Let $E$
be a
locally trivial QVB over $B$ corresponding to this covering.\\
Then $E$ is the space of 
sections of a locally trivial vector bundle $N$. \end{pr}
Proof: We want to give here only the idea of the proof. To prove the first 
assertion one defines the module homomorphisms $\zeta_i$ in terms of the 
trivializations $\psi_i: N \longrightarrow U_i \times V$ by \[ \zeta_i(s) =
\psi_i \circ s \] (identifying $\psi_i(s(x))=(x,v)$ with $v$) and shows the
conditions claimed for $\zeta_i$.
\\To prove the second assertion one construct in term of the given locally 
trivial QVB a locally trivial vectorbundle in the following way. Let 
$X=\bigcup_i \{i\} \times U_i \times V$. One obtains a locally trivial
vector 
bundle 
$N$ over $M$ by factorizing $X$ with respect to the following relation: 
 \[ (i,x,v) \sim (j,x',v')~~\mbox{if}~~ x=x'~~\mbox{and}~~ 
 v= \phi_{{ij}_E}(1 \otimes v')(x). \] 
 One proves that the module of sections $\Gamma(N)$ is isomorphic to $E$. 
 \hfill$\Boxneu$. 
 \\\\ Assume that there is given an LC differential algebra $\Gamma(B)$
with complete covering $(J_{i_{\Gamma}})_{i \in I}$ such that 
 $pr_0(J_{i_{\Gamma}})	=J_i$, i.e. the 
 factor algebras $\Gamma(B)/J_{i_{\Gamma}}=\Gamma(B_i)$ are differential 
 calculi over the  factor algebras $B_i=B/J_i$. One can construct a locally 
 trivial QVB $((E_{\Gamma},\Gamma(B),\kappa_{\Gamma}),V,(\zeta_{i_{\Gamma}},
 J_{i_{\Gamma}})_{i \in I} )$ from $E$ in the following way. 
One extends $\phi_{{ij}_{E}}$ to $\Gamma(B_{ij}) 
\otimes V$ by \[ \phi_{{ij}_{E}}(\gamma \otimes v)=
\gamma\phi_{{ij}_{E}}(1 \otimes v),\;\; \gamma \in
\Gamma(B_{ij}),\; 
v \in V.\] In terms of this extended module isomorphism one defines 
$E_{\Gamma}$ by \begin{equation} E_{\Gamma}:=\{ (e_{i_{\Gamma}})_{i \in I}
\in 
\bigoplus_{i \in I} \Gamma(B_i) \otimes V| (\pi^i_{j_{\Gamma}} \otimes 
id)(e_{i_{\Gamma}})=\phi_{{ij}_{E}} \circ (\pi^j_{i_{\Gamma}} \otimes 
id)(e_{i_{\Gamma}}) \}, \end{equation}
where the homomorphisms 
$\pi^i_{j_{\Gamma}} : \Gamma(B_i) \longrightarrow 
\Gamma(B_{ij}):=\Gamma(B)/(J_{i_{\Gamma}} 
+ J_{j_{\Gamma}})$ are the canonical projections.
\\By defining $\kappa_{\Gamma} : \Gamma(B) \longrightarrow End(E_{\Gamma})$
as 
\[ \kappa_{\Gamma}(\gamma)((e_{i_{\Gamma}})_{i \in 
I}):=(\pi_{i_{\Gamma}}(\gamma)e_{i_{\Gamma}})_{i \in I} \]
 ($\pi_{i_{\Gamma}} : \Gamma(B) \longrightarrow 
  \Gamma(B_i)=\Gamma(B)/J_{i_{\Gamma}}$ is 
  the canonical projection.) and $\zeta_{i_{\Gamma}} : E_{\Gamma} 
  \longrightarrow 
  \Gamma(B_i) \otimes V$ as the i-th projection one obtains a locally
trivial QVB   $((E_{\Gamma}, \Gamma(B), \kappa_{\Gamma}),V,( 
  \zeta_{i_{\Gamma}},J_{i_{\Gamma}})_{i \in I} )$. \\ 
  Now one can consider connections on such locally trivial QVB.
We add to the usual definition of a connection in a ``vector bundle'' as
a covariant derivative (\cite{co8.5}) a condition of compatibility with the
covering. \begin{de} Let 
  $(E_{\Gamma},\Gamma(B),\kappa_{\Gamma}),V,(\zeta_{i_{\Gamma}},J_i)_{ i
\in   I} )$ be the locally trivial QVB just defined. 
  A connection on $E_{\Gamma}$ is a linear map $\nabla: E_{\Gamma} 
  \longrightarrow E_{\Gamma}$ satisfying 
  \begin{eqnarray} \label{nab}	\nabla(\gamma e) &=& (d\gamma) e + (-1)^n
\gamma \nabla(e),\; 
  \gamma \in \Gamma^n(B),\; e \in E_{\Gamma} \\  \label{nabi}
  \nabla (ker \zeta_{i_{\Gamma}}) & \subset & ker \zeta_{i_{\Gamma}},\;
\forall i \in I. \end{eqnarray} \end{de}
  It is easy to see that from the property (\ref{nabi}) follows that there
exist 
  connections $\nabla_i$ on $\Gamma(B_i) \otimes V$ and $\nabla^i_{ij}$ on 
  $\Gamma(B_{ij}) \otimes V$ such that \begin{eqnarray*} \nabla_i \circ 
  \zeta_{i_{\Gamma}} &=& \zeta_{i_{\Gamma}} \circ \nabla \\
  \nabla^i_{ij} \circ (\pi^i_{j_{\Gamma}} \otimes id) &=&
(\pi^i_{j_{\Gamma}} \otimes id) \circ \nabla_i. \end{eqnarray*}
  \begin{pr} \label{conne} Connections on $E_{\Gamma}$ are in one to one 
  correspondence with 
  families of connections $\nabla_i : \Gamma(B_i) \otimes V \longrightarrow
  \Gamma(B_i) \otimes V$ satisfying 
  \begin{equation} \label{nabtriv} \nabla^i_{ij} = \phi_{{ij}_E} \circ 
  \nabla^j_{ij} \circ \phi_{{ji}_E}. \end{equation} \end{pr}
  Proof: Let $\nabla$ be a connection on $E_{\Gamma}$. There exists a
family of 
  connections $\nabla_i$ on $\Gamma(B_i) \otimes V$ such that \[ \nabla 
  ((e_{i_{\Gamma}})_{i \in I})=(\nabla_i(e_{i_{\Gamma}}))_{i \in I}. \]
  Because of the identities \begin{eqnarray*}
  (\pi^i_{j_{\Gamma}} \otimes id)(e_{i_{\Gamma}})&=& \phi_{{ij}_E} \circ 
  (\pi^j_{i_{\Gamma}} \otimes id)(e_{j_{\Gamma}})\\
  (\pi^i_{j_{\Gamma}} \otimes id) \circ \nabla_i(e_{i_{\Gamma}})&=& 
  \phi_{{ij}_E} \circ (\pi^j_{i_{\Gamma}} \otimes id) \circ
\nabla_j(e_{j_{\Gamma}}) 
  \end{eqnarray*} one obtains \begin{eqnarray*}
  \nabla^i_{ij} \circ (\pi^i_{j_{\Gamma}} \otimes id)(e_{i_{\Gamma}}) &=& 
  \phi_{{ij}_E} \circ \nabla^j_{ij} \circ (\pi^j_{i_{\Gamma}} \otimes 
  id)(e_{j_{\Gamma}}) \\
  &=& \phi_{{ij}_E} \circ \nabla^j_{ij} \circ \phi_{{ji}_E} \circ 
  (\pi^i_{j_{\Gamma}} 
  \otimes id)(e_{i_{\Gamma}}). \end{eqnarray*}
  Conversely, if there is given a family of connections $\nabla_i$
satisfying 
  property (\ref{nabtriv}) the image of the linear map $\nabla$ defined by 
  \[ \nabla((e_{i_{\Gamma}})_{i \in I})=(\nabla_i(e_{i_{\Gamma}}))_{i \in
I} 
  \] lies in $E_{\Gamma}$ and has the properties of a connection.
\hfill$\Boxneu$
  \\Remark: Let the family $(e_i)_{i \in I}$ be a linear basis in $V$. Let $\nabla$
be a 
  connection on $\Gamma(B) \otimes H$. Then there is a 
  family of one forms $(A^j_i)$ such that $\nabla(1 \otimes e_i)=\sum_j
A^j_i
  \otimes e_j$.
  \begin{de} Let $\nabla$ be a connection on $E_{\Gamma}$. We call the
linear 
  map 
  $\nabla^2$ the curvature of $\nabla$. \end{de}
  Note that $\nabla^2(\gamma e)=\gamma \nabla^2(e)$.
  \\In the sequel we will be interested in QVB associated to a QPFB.
  We define these as follows (see also \cite{Kon}):
  \begin{de} Let $\cal P$ be a locally trivial QPFB, let $F$ be a vector space 
   and let $\rho: F  \longrightarrow H \otimes F$ be a left $H$ coaction.
  The vector bundle $E({\cal P},F)$ associated to $\cal P$ and $\rho$ is
defined as   \begin{equation} E({\cal P},F):= \{ e \in {\cal P} \otimes F|
(\Delta_{\cal P} 
  \otimes id)(e)=(id \otimes \rho)(e) \}. \end{equation} \end{de}
  Remark: This is also called co-tensor product of $\cal P$ and $F$.
  \begin{pr} The associated vector bundle $E({\cal P},F)$ is a locally
trivial QVB.  \end{pr}
  Proof: By formula (10) of \cite{cama1}. $E({\cal P},F)$ has the form
\begin{eqnarray*} E({\cal P},F) = \{ ( e_i)_{i \in I} \in 
  \bigoplus_{i \in I} B_i \otimes H 
  \otimes F| (\pi^i_j \otimes id \otimes id)(e_i) &=& (\phi_{ij} \otimes
id) \circ   (\pi^j_i \otimes id \otimes id)(e_j); \\ (id \otimes \Delta \otimes 
  id)(e_i)&=&(id   \otimes id \otimes \rho)(e_i) \}. \end{eqnarray*}
  There are isomorphisms \[ id \otimes \varepsilon \otimes id: E_i:=\{
e_i \in B_i 
  \otimes H \otimes F| id \otimes \Delta \otimes id(e_i) = id
\otimes id   \otimes \rho(e_i) \} \longrightarrow B_i \otimes F, \]
  where $(id \otimes \varepsilon \otimes id)^{-1} = id \otimes \rho$,
such that 
  all $ e_i \in E_i$ are of the form \[\sum \sum_k a_k \otimes f_{k_{(-1)}}
  \otimes f_{k_{(0)}}.\]
  One easily verifies that 
$((\tilde{E},B,\tilde{\kappa}),F,(\tilde{\zeta_i},
  J_i)_{i \in I})$ defined by 
  \[ \tilde{E}({\cal P},F) := \{ ( \sum_k a^i_{k_i} \otimes f^i_{k_i} \}
\in   \bigoplus_{i \in I} B_i \otimes F|\sum_k \pi^i_j(a^i_{k_i}) \otimes
f^i_{k_i}= \sum \sum_k \pi^j_i(a^j_{k_i}) \tau_{ji}(f^j_{k_{j(-1)}} )
\otimes 
  f^j_{k_{j(0)}} \} \]	 \begin{eqnarray*} \kappa(a)(( \tilde{e}_i )_{i \in
I}) &=& (\pi_i(a) 
  \tilde{e}_i )_{i \in	 I},\; a \in B,\; ( \tilde{e}_i )_{i \in I} \in 
  \tilde{E} \\
  \zeta_i && \mbox{i-th projection}  \end{eqnarray*}  is a locally trivial
QVB (For the definition of the transition functions $\tau_{ij}$ see 
\cite{cama1}). In terms of the isomorphisms $id \otimes \varepsilon 
\otimes id$ one
can define a module isomorphism $\epsilon: E({\cal P},F) \longrightarrow 
  \tilde{E}({\cal P},F)$ by \[ \epsilon((e_i)_{i \in I}):= (id \otimes \varepsilon \otimes
  id(e_i))_{i \in I}. \]
  This isomorphism exists due to the glueing properties.
  This is easy to see: An element  $(e_i)_{i \in I} \in E({\cal P},F)$ 
  is of the 
  form $(e_i)_{i \in I}=(\sum \sum_{k_i} a^i_{k_i} \otimes f^i_{k_{i(-1)}} 
  \otimes f^i_{k_{i(0)}})_{i \in I}$
  satisfying \[ \sum \sum_{k_i}\pi^i_j(a^i_{k_i}) \otimes f^i_{k_{i(-1)}} 
  \otimes f^i_{k_{i(0)}} = 
  \sum \sum_{k_j} \sum_{k_j} \pi^j_i(a^j_{k_j})\tau_{ji}(f^j_{k_{j(-2}})
\otimes 
  f^j_{k_{j(-1)}} \otimes 
  f^j_{k_{j(0)}}.\] (see formula (10) of \cite{cama1})
  \\Applying $id \otimes \varepsilon \otimes id$ to both sides of 
  this equation, one obtains for the element $(id \otimes \varepsilon 
  \otimes 
  id(e_i))_{i \in I}$ the property 
  \[ \sum_{k_i} \pi^i_j(a^i_{k_i}) \otimes f^i_{k_i} = \sum \sum_{k_j} 
  \pi^j_i(a^j_{k_j})\tau_{ji}(f^j_{k_{j(-1)}}) \otimes f^j_{k_{j(0)}},\]
  i.e. the image of $\epsilon$ lies in $\tilde{E}$. The inverse of
$\epsilon$ is 
  the map $(\tilde{e}_i)_{i \in I} \longrightarrow (id \otimes 
  \rho(\tilde{e}_i))_{i \in I}$, i.e $\epsilon$ is an isomorphism.
  We obtain the locally trivial QVB 
  $((E({\cal P},F),B,\kappa),F,(\zeta_i,J_i)_{i \in I})$ 
  with	\begin{eqnarray*} \kappa &=& \tilde{\kappa} \circ \epsilon \\
  \zeta_i &=& \tilde{\zeta}_i \circ \epsilon.\;\;\; \end{eqnarray*} 
  {}\hfill$\Boxneu$  
  \\\\
As in the general case of locally trivial vector 
  bundles   one can construct a locally trivial vector bundle 
  $((E_{\Gamma}({\cal 
P},F),\Gamma_m(B),\kappa_{\Gamma}),(\zeta_{i_{\Gamma}},ker 
 \pi_{i_{\Gamma_m}} )_{i \in I}) $ from the associated vector bundle 
  $E({\cal P},F)$. The LC differential algebra $\Gamma_m(B)$ is the maximal 
  embeddable LC-differential algebra induced
from the differential structure $\Gamma_c({\cal P})$ on the locally trivial QPFB
${\cal P}$ (see \cite{cama} and \cite{cama1}).
  Let $\phi_{{ij}_E}: B_{ij} \otimes F \longrightarrow B_{ij} \otimes F$ be
  defined by 
\[ \phi_{{ij}_E}(a \otimes f):=a \tau_{ji}(f_{(-1)}) \otimes 
  f_{(0)}.
\]   
By definition, 
\begin{equation} 
E_{\Gamma}({\cal P},F)=
  ( (e_{i_{\Gamma}})_{i \in I} \in \bigoplus_{i \in I} \Gamma(B_i) \otimes
F| (\pi^i_{j_{\Gamma_m}} \otimes id)(e_{i_{\Gamma}})=\phi_{{ij}_{E}} \circ 
  (\pi^j_{i_{\Gamma_m}} \otimes id)(e_{j_{\Gamma}}) \}. 
\end{equation}

\begin{pr} 
Let 
\[
 hor E({\cal P}, F):=\{ \gamma_E \in hor \Gamma_c({\cal
P}) \otimes F| (\Delta_{{\cal P}_{\Gamma}} \otimes id)(\gamma_E) = (id
\otimes \rho)(\gamma_E) \}. 
\]
  $hor E({\cal P},F)$ and $E_{\Gamma}({\cal P},F)$ are isomorphic as 
left $\Gamma_m(B)$-modules. 
\end{pr}
  Remark: $hor E({\cal P},F)$ is in the classical situation the space of 
  horizontal forms ``of type $\rho$''.
  \\\\Proof: By definition of $hor \Gamma_c({\cal P})$, 
\begin{eqnarray*} hor E({\cal P},F)&=& \{ (\gamma_{i_E})_{i \in I} \in 
  \bigoplus_{i \in I} \Gamma(B_i) \hat{\otimes} H \otimes F| \\
 && ((\pi^i_j \otimes id)_{\Gamma} \otimes id)(\gamma_{i_E}) = 
  (\phi_{{ij}_{\Gamma}} \otimes id) \circ ((\pi^j_i 
  \otimes id)_{\Gamma} \otimes	 id)(\gamma_{j_E}); \\
&&(id \otimes \Delta \otimes id)(\gamma_{i_E}) = 
  (id \otimes id \otimes \rho)(\gamma_{i_E}) \}. \end{eqnarray*}
  The last condition means that the i-th components have the 
  form \[ \gamma_{i_E}= \sum \sum_{k_i} \gamma^i_{k_i} \otimes
f^i_{k_{i(-1)}} \otimes 
  f^i_{k_{i(0)}}. \] Now one defines again an isomorphism
$\epsilon_{\Gamma} : 
  hor E({\cal P},F) \longrightarrow E_{\Gamma}({\cal P},F) $ by 
  \[  \epsilon_{\Gamma}((\gamma_{i_E})_{i \in I})=(id \otimes \varepsilon 
  \otimes id(\gamma_{i_E}) )_{i \in I}. \] 
This isomorphism exists, 
  if one can show that from the gluing conditions in $hor E({\cal P},F)$
the gluing 
  conditions in $E_{\Gamma}({\cal P},F)$ follow. To this end apply $P_{inv}
  \otimes id$ (formula (66) of \cite{cama1}) to
  \[ \sum \sum_{k_i} (\pi^i_j\otimes id)_{\Gamma}(\gamma^i_{k_i} \otimes 
  f^i_{k_{i(-1)}}) \otimes 
  f^i_{k_{i(0)}}=\sum \sum_{k_j} \phi_{{ij}_{\Gamma}} \circ 
  (\pi^j_i \otimes id)_{\Gamma}(\gamma^j_{k_j} \otimes f^j_{k_{j(-1)}})
\otimes 
  f^j_{k_{j(0)}}. \]  By the definition
$\pi^i_{j_{\Gamma_g}}(\gamma)=(\pi^i_j 
  \otimes id)_{\Gamma}(\gamma \hat{\otimes} 1)$ one 
  obtains the gluing condition in $E_{\Gamma}({\cal P},F)$, i.e. 
  \[ \sum_{k_i} \pi^i_{j_{\Gamma_g}}(\gamma^i_{k_i}) \otimes f^i_{k_i} = 
  \sum_{k_j} 
  \phi_{{ij}_{E_{\Gamma}}}(\pi^j_{i_{\Gamma_g}}(\gamma^j_{k_j}) \otimes 
  f^j_{k_j}),\] which means that the image of $\epsilon_{\Gamma}$ lies in 
  $E_{\Gamma}({\cal 
  P},F)$. Conversely, the inverse of $\epsilon_{\Gamma}$ is obviously
defined by 
  \[ \epsilon^{-1}_{\Gamma}((e_{i_{\Gamma}} )_{i \in I})= (id \otimes 
  \rho(e_{i_{\Gamma}}))_{i \in I}. \]	 \hfill$\Boxneu$
  \\\\Now we are interested in connections on such associated vector
bundles. 
  An important class of connections are the connections induced from left 
  left covariant derivations on the locally trivial QPFB.
  
\begin{pr} Every left covariant derivation on the locally trivial QPFB 
  $\cal P$ determines 
  uniquely a connection on $E_{\Gamma}({\cal P},F)$. 
\end{pr}
  Proof: One defines a linear map $\nabla: E_{\Gamma} \longrightarrow 
  E_{\Gamma}$ by 
\[\nabla:=\epsilon_\Gamma \circ (D_l \otimes id) \circ \epsilon_\Gamma^{-1},\]
which is easily seen to be a connection on 
  $E_{\Gamma}({\cal P},F)$. \hfill$\Boxneu$
  \\\\
Remark: Because of the bijection between left and right covariant 
  derivations also right covariant derivations on $\cal P$ determine 
  connections on $E_{\Gamma}({\cal P},F)$.
  \\
The curvature of such a connection is related to the 
  curvature on the locally trivial QPFB $\cal P$ by  
  \[\nabla^2= \epsilon_{\Gamma} \circ (D_l^2 \otimes id) \circ 
  \epsilon^{-1}_{\Gamma}. \]
The corresponding connections  
  $\nabla_i: \Gamma(B_i) \otimes F \longrightarrow \Gamma(B_i) \otimes F$
have the form 
\[
 \nabla_i(\gamma \otimes f)= d\gamma \otimes f - \sum 
  (-1)^n \gamma A_{l_i}(f_{-1}) \otimes f_0, \; \gamma \in
\Gamma^n(B_i),\;f \in F. 
\] 
The curvatures of these connections are 
  \[
\nabla^2_i(\gamma \otimes f)=\sum \gamma F_{l_i}(f_{-1}) 
  \otimes f_0.
\]


\begin{thebibliography}{100}
\bibitem{brma} Brzezi\'nski, T., and S. Majid: Quantum group gauge theory
on quantum spaces, {\it Commun. Math. Phys.} {\bf 157} (1993), 591-638,
{\sf hep-th/9208007}, 
Preprint DAMTP/92-27
\bibitem{Kon} Budzy\'nski, R. J. and W. Kondracki: Quantum principal fiber 
bundles: Topological aspects, {\it Rep. Math. Phys.} {\bf 37} (1996),
365-385, 
 preprint 517 PAN Warsaw 1993, {\sf hep-th/9401019}
\bibitem{cama} Calow, D. and R. Matthes: Covering and gluing of algebras
and differential algebras, 
{\it J. Geom. and Phys.} {\bf 32} (2000), 364--396,
{\sf math.QA/9910031}, Preprint
NTZ 25/1998
\bibitem{cama1} Calow, D. and R. Matthes: Connections on locally trivial
quantum principal fibre bundles, submitted to {\it J. Geom. and Phys.},
{\sf math.QA/0002228}
\bibitem{co8.5} Connes, A.: {\it Noncommutative Geometry}, Academic Press,  
Inc., San Diego, New York etc., 1994
\bibitem{dur1} Durdevic, M.: Geometry of quantum principal bundles II,
{\it Rev. Math. Phys.} {\bf 9} (5) (1997), 531--607, {\sf q-alg/9412005}
\bibitem{dur1.1} Durdevic, M.: Quantum principal bundles and Tannaka-Krein 
duality theory, {\it Rep. Math. Phys.} {\bf 38} (1996), 313--324
\bibitem{hama1} Hajac, P. M. and R. Matthes: Frame, cotangent and tangent
bundles of the quantum plane, in Proc. II Int. Workshop
{\it Lie Theory and its Applications in Physics II}, ed. by H.-D. Doebner, V.
K. Dobrev \& J. Hilgert, World Scientific, Singapore 1998, 377--385,
{\sf math.QA/9803127}, DAMTP-98-17, NTZ 9/98
\bibitem{pf1} Pflaum, M. J.: Quantum groups on fibre bundles,
{\it Commun. Math. Phys.} {\bf 166} (1994), 279-315, {\sf hep-th/9401085}
\end{thebibliography}
\end{document}